\documentclass [12pt]{article}
\usepackage{xcolor}
\usepackage{amsfonts}
\usepackage{amssymb,amsmath,amsthm}
\usepackage{graphicx}
\usepackage{multirow,booktabs}

\usepackage[authoryear]{natbib}
\newtheorem{theorem}{Theorem}[section]

\newtheorem{corollary}{Corollary}[section]

\newtheorem{example}{Example}[section]

\newcommand{\dint}{\displaystyle\int}
\newcommand{\diag}{\mbox{diag}}

\newcommand{\crno}{\cr\noalign{\vskip1.5mm}}
\def \bdm#1{\mbox{\boldmath$#1$}}

\def \bLam {\mbox{\boldmath $\Lambda$}}
\def \btheta {\mbox{\boldmath $\theta$}}
\def\c{{\boldsymbol c}}

\def\f{{\boldsymbol f}}
\def\g{{\boldsymbol g}}
\def\h{{\boldsymbol h}}
\def\x{{\boldsymbol x}}

\def\z{{\boldsymbol z}}
\def\C{{\boldsymbol C}}
\def\D{{\boldsymbol D}}

\def\M{{\boldsymbol M}}
\def\N{{\boldsymbol N}}

\def\Q{{\boldsymbol Q}}

\def\Z{{\boldsymbol Z}}

\newcommand{\trace}{\mbox{trace}}

\begin{document}

\title{ Design admissibility and de la Garza phenomenon in multi-factor experiments}

\author{
{\small Holger Dette} \\
{\small Ruhr-Universit\"at Bochum} \\
{\small Fakult\"at f\"ur Mathematik}\\
{\small Bochum, Germany} \\
\and
{\small Xin Liu } \\
{\small  College of Science} \\
{\small  Donghua University}\\
{\small  Shanghai 201600, China } \\
\and
{\small Rong-Xian Yue} \\
{\small Department of Mathematics,} \\
{\small   Shanghai Normal University}\\
{\small  Shanghai 200234, China}\\
}

  \maketitle

\begin{abstract}

The determination of an optimal design for a given regression problem  is an intricate optimization problem, especially
for models with multivariate predictors.
% How to reduce the complexity of the problem is the key to the successful solution.
Design admissibility and invariance are main tools to  reduce the complexity of the  optimization problem
 and have  been successfully applied for models with univariate predictors.
 In particular several authors have developed sufficient conditions for the existence of saturated designs in univariate models, where the number of
 support  points of the optimal design equals the number of parameters.
These  results generalize the celebrated de la Garza phenomenon  \citep{Garza1954} which states that for a polynomial regression model of degree $p-1$ any optimal design can be based on at most
$p$  points.

This paper provides - for the first time -    extensions of these results  for models with a multivariate predictor.
In particular  we study  a geometric characterization of the support points of an optimal design   to  provide sufficient conditions for the occurrence of the de la Garza phenomenon in models with  multivariate predictors
and   characterize properties of admissible designs in terms of admissibility of designs in {\it conditional} univariate regression models.
\end{abstract}

\vspace{9pt}
\par

\vspace{9pt}

\section{Introduction}
 \label{sec1}
\def\theequation{1.\arabic{equation}}
\setcounter{equation}{0}

It is well known that an appropriate choice of an experimental design can improve the quality of statistical analysis
substantially, and therefore the problem of constructing optimal designs for regression models has found considerable attention in
the literature \citep[see, for example, the monographs of][]{Puke1993,Randall2007}. However, the determination of
 an optimal design often results in an intricate optimization problem that is difficult to handle, in particular for models used for experiments with multivariate predictors.

 A useful strategy is to simplify the  problem by identifying  subclasses of  relatively simple designs, which must contain
 the optimal design. A prominent example of such a class is the class of admissible designs consisting of the designs with an information matrix,
 that cannot be improved  by an information matrix of another design with respect to the Loewner ordering.  In  decision theoretic terms the set of admissible designs
 therefore   forms a   {\it complete class}, in the sense that the information matrix of any  inadmissible design  may be improved by the
 information matrix of  an admissible design.  It is well known that  optimal
 designs with respect to   the most of  the commonly used optimality criteria  must
  be admissible  \citep[see][Chapter 10.10]{Puke1993} and consequently  in these cases
  the determination of optimal designs  can be restricted to the  class of admissible designs.
  Along this line, in a series of remarkable papers  \cite{YS2009,YS2012}, \cite{Yang2010}, \cite{DM2011}, \cite{DS2013}  and \cite{hu2015}
  derived several complete classes of designs for regression models with a univariate predictor. In particular it is demonstrated that the
  celebrated de la Garza phenomenon \citep{Garza1954}, which states that for a polynomial regression model of degree $p-1$ any optimal design can be based on at most
$p$  points,   appears in a broad class of regression models with a univariate predictor.

While these  methods  provide   a very powerful tool for the determination of optimal designs, its  application is limited to single-factor experiments
since the key tools to prove these results are not available for functions of several variables. For example, the characterizations developed in \cite{DM2011} and  \cite{DS2013}
are based the theory of Chebyshev systems \citep[see][]{karstu1966}, which  requires regression functions with a univariate argument.
Consequently, for  regression models with   multivariate predictor optimal design problems, including investigations of admissibility, have been  mostly treated on a case-by-case
analysis  using various techniques.  For example,  \cite{Heiligers1992}  investigated admissible experimental designs in a multiple polynomial regression model.
{\cite{YZH2011} derived a class of admissible  designs for the commonly used multi-factor logistic and probit models.
\cite{HHL2020} characterized an essentially complete class with respect to Schur ordering for binary response models with multiple nonnegative explanatory variables.}
Moreover, for several specific models with a multivariate predictor  optimal designs with respect to various criteria have been determined.
Exemplarily, we mention   \cite{GGHS2007}, who  studied locally $D$-optimal designs for  generalized linear models using a canonical transformation,
  \cite{BDW2011}, who  showed that  in additive partially nonlinear models
  $D$-optimal designs can be found as the products of the corresponding $D$-optimal designs in one dimension, \cite{dette2014}, who
studied $E$-optimal designs for second order response surface models,
    % \cite{GMS2017,GMS2018,GMS2019}
 \cite{GMS2018}, who  discussed  locally $D$-optimal designs for the Cobb-Douglas model, \cite{KHN2018}, who investigated $D$-optimal designs for the two-variable binary logistic regression model with interaction, and \cite{CGHHL2019}, who used the moment-sum-of-squares hierarchy
of semidefinite programming problems to solve approximate optimal design problems for multivariate polynomial regression on a compact space.

In the present  paper we study these  problems
from a more general point of view. In particular  we develop  a
geometric characterization for  the support points of an optimal design which can be used to  derive sufficient conditions for the  occurrence of  the
de la Garza phenomenon in regression models with a multivariate predictor.  Our strategy is to handle the design problem by considering  the  dual
optimization problem. Moreover, in  contrast to the previous literature, which considers  characterizations in terms of the
explanatory variable, our approach uses the {\it induced design space} of the regression model under consideration.
Moreover, we also provide a necessary  condition for a class of designs to be admissible
in terms of the admissibility   of  the designs in the  corresponding {\it conditional models}.

In Section \ref{sec2} we develop  sufficient conditions for the occurrence of the de la Garza phenomenon based on the geometric
 characterization of the support points of an optimal design. Section \ref{sec3} introduces the concept of conditional models and   designs,
 which are used to investigate design admissibility for models with multivariate predictors. In Section 4, we illustrate the  potential of our approach
 in three examples considering various nonlinear models with a multivariate predictor. Finally all proofs of our technical results are deferred to Section \ref{sec5}.

\section{Optimal designs  and  a geometric characterization}
 \label{sec2}
\def\theequation{2.\arabic{equation}}
\setcounter{equation}{0}

We begin stating the optimal design problem as considered, for example,  in  \cite{Puke1993}.
Throughout this paper let  Sym($s$)  denote  the set of all real symmetric  $s \times s$ matrices, NND($s$) $\subset $  Sym($s$)   the set of all nonnegative definite matrices
and  PD($s$) $\subset $  NND($s$) the set of positive definite matrices.
We consider the common linear regression model
\begin{eqnarray}\label{LM}
y= \f^\top(\x)\btheta + \varepsilon,
\end{eqnarray}
where  $\x=(x_1,\cdots,x_q)^{\top}$ is a $q$-dimensional vector of  predictors which varies in the design space $\mathcal{X}\subset\mathbb{R}^q$,  $\f(\x)$ is a $k$-dimensional  vector of known linearly independent
regression functions,  $\btheta \in \mathbb{R}^{k}$ denotes the   vector of unknown parameters, and  $\varepsilon$ is a random
variable with mean $0$ and constant variance $\sigma^2 >0$. We assume that the experimenter can take $n$ independent observations of the form
$y_{i}= \f^\top(\x_{i})\btheta + \varepsilon_{i}$ ($i=1, \ldots , n$) at experimental conditions $\x_{1}, \ldots ,  \x_{n}$.

Following \cite{Kiefer1974}  we define   a (approximate) design for model (\ref{LM})  as a probability measure $\xi$ on the design space $\mathcal{X}$ with finite support and the information matrix of the design  $\xi$ in model \eqref{LM}
 by
 \begin{equation}\label{Informatix Matrix}
{\M}(\xi)=\dint_{\mathcal{X}}\f({\x}){\f}^\top({\x})\xi(d\x).
\end{equation}
If the design $\xi$ has  masses $\xi_{1}, \ldots , \xi_{m}$ at   $m$ support points $\x_1,\ldots , \x_m$, and $n$ observations can be taken, the quantities $\xi_{\ell} n$ are rounded to non-negative
integers, say $n_{\ell}$, such that $\sum_{\ell=1}^{m} n_{\ell} =n$ and the experimenter takes $n_{\ell}$ observations at each $\x_{\ell}$ ($\ell =1  \ldots , m$). In this case the covariance matrix of the least squares estimator
$ \sqrt{n} \hat \btheta $  for the parameter $ \btheta $ in model \eqref{LM}
converges  to the matrix $\sigma^{2}  {\M}^{{-1}}(\xi)$. which is used to measure the accuracy of the estimator $ \hat \btheta $.

We use the notation $\Xi$ for the set of all approximate designs on the design space ${\cal X}$  and $\mathcal{M}(\Xi) = \{ \M (\xi ) ~|~\xi  \in \Xi \} $ for the set of all information matrices.
 An optimal design $\xi^{*}$ maximizes an appropriate function, say $\phi$,  of the information matrix $M(\xi)$, where  $\phi$ : NND($s$)$\rightarrow\mathbb{R}$
 is a positively homogeneous, super-additive, nonnegative, non-constant and upper semi-continuous function.
Throughout this paper we   call a function  with these properties {\it optimality criterion} or  {\it information function}.
The most prominent optimality criteria are the matrix means defined by
\begin{equation}\label{phi_q}
\phi_p(\C)= \left\{\begin{array}{ll}
  \big (\frac{1}{s}\trace(\C^p)\big)^{1/p}   & \textmd{for}\; p\in(-\infty, 1] \setminus\{0\} \crno
   (\det(\C))^{1/s}                                            & \textmd{for}\;p=0 \crno
   \lambda_{\min}(\C)                                    & \textmd{for}\;p=-\infty
 \end{array},
 \right.
\end{equation}
which include the classical $A$-, $D$- and $E$-optimality criteria as special cases $p=-1$, $p=0$ and $p=-\infty$, respectively
(here we define $\phi_p(\C) =0 $ if  $\C\in$ NND($s$)$\setminus$PD($s$)).

Given an optimality criterion $\phi$ on NND($k$) the {\it design problem} then reads as follows
\begin{eqnarray}\label{DEP}
\max_{\M\in\mathcal{M}(\Xi)} \phi(\M),
\end{eqnarray}
where, in a second step, one has to identify a design $\xi^{*}$ corresponding to a maximizer $\M^{*}$ of \eqref{DEP}. Any design with this property is called {\it $\phi$-optimal design}. As pointed out in the introduction, an important problem
in optimal design theory is to identify  sufficient conditions on the regression model \eqref{LM}  such that (approximate) optimal    designs are saturated, which means that the number of support points of the design
coincides with the dimension of the parameter. This property   is called de la Garza phenomenon referring to the  famous result of \cite{Garza1954}, which shows that the $G$-optimal design in a polynomial regression of degree $p-1$ on a compact interval  
has $p$ support points.  While this problem has found considerable attention for models with one-dimensional predictors (see the references mentioned in the introduction), there are
- to our best knowledge -  no general results available which characterize saturated designs in   models with a multivariate predictor.

We begin with a  geometric characterization of the support points of a $\phi$-optimal design, which can be used to derive  sufficient conditions for the occurrence of the de la Garza phenomenon in models with multivariate predictors.
For this purpose we define for a matrix $\Z  =(\z_{1},\ldots , \z_{k} ) \in \mathbb{R}^{k \times k}$ a linear transformation $\h_{\Z} :   \mathbb{R}^{k} \to  \mathbb{R}^{k}$ by
\begin{equation} \label{hd3}
\h_{\Z}(\x):=(h_{\Z 1}(\x), \ldots, h_{\Z k}(\x))^\top :=  {\Z}^\top\f(\x)=({\z_1}^\top\f(\x), \ldots, {\z_k}^\top\f(\x))^\top
\end{equation}
and consider the corresponding point
\begin{equation} \label{hd3}
P_{\Z} (\x) =(h_{\Z 1}^2(\x), \ldots, h_{\Z k}^2(\x))^{\top} \in \mathbb{R}^{k}~.
\end{equation}

\begin{theorem}\label{thm1}
Let $\xi^*=\{(\x^*_i, w^*_i)\}_{i=1}^n$ be a $\phi$-optimal design for the  regression model (\ref{LM}). There exists an orthogonal matrix, say $\Z^*=(\z^*_{1}, \dots, \z^*_{k}) \in \mathbb{R}^{k \times k}$
 with a linear transformation $\h_{\Z^*}$ of the form \eqref{hd3}, such that the vectors  {$P_{\Z^{*} } (\x_{1}^*), \ldots , P_{\Z^{*}} (\x_{n}^*)$}
define at most $k$ different supporting hyperplanes of the  $k$-dimensional polytope
\begin{eqnarray}\label{hd4}
{\cal P}_{\Z^*} :=
\big \{ \lambda= (\lambda_1, \ldots, \lambda_k)^{\top} \colon \lambda_i\geq0,  \;\; \forall   i=1, \ldots, k, \;P_{\Z^*}^{\top} (\x) \lambda
% = \h_{\Z 1}^2(\x)\lambda_1+\cdots+\h_{\Z k}^2(\x)\lambda_k
\leq 1 \;\; \forall \x\in\mathcal{X} \big \}.
\end{eqnarray}
Moreover, if  $\f(\x^*_{i})$  and $\f(\x^*_{j})$ are two vectors corresponding to the same supporting hyperplane they have the same length.
\end{theorem}

\begin{example}  \label{ex1}
{\rm  To illustrate the result given by Theorem $\ref{thm1}$, we  consider a linear regression in two variables with no intercept, that is
$\f(\x)=(x_1, x_2)^\top$, where $\x=(x_1, x_2)\in[0, 1]^2$. If
$$
\Z=\left(\begin{matrix}
  \cos t&\sin t\\
  -\sin t&\cos t
\end{matrix}\right)
$$
is a  $2 \times 2$  orthogonal matrix, then  the vector $\h_\Z$  in \eqref{hd3} is given by
$$
\h_\Z(\x)=(x_1\cos t-x_2\sin t, x_1\sin t+x_2\cos t)^\top,
$$
and it is  easy to see that the polytope
\begin{equation} \label{hd5}
{\cal P}_{\Z}= \big\{(\lambda_1, \lambda_2): \lambda_1\geq0, \lambda_2\geq0, (x_1\cos t-x_2\sin t)^2\lambda_1+(x_1\sin t+x_2\cos t)^2\lambda_2\leq 1 \;\; \forall \x\in\mathcal{X}\big\}
\end{equation}
is determined by at most three half planes, which are defined by
\[
\begin{array}{rcr}
\lambda_1\cos^2 t+\lambda_2\sin^2 t&\leq& 1,\crno
\lambda_1\sin^2 t+\lambda_2\cos^2 t&\leq& 1,\crno
\lambda_1(1-\sin 2t)+\lambda_2(1+\sin 2t)&\leq& 1,
\end{array}
\]
and correspond  to the points $(1, 0), (0, 1)$ and $(1, 1)$, respectively.
 Therefore, the support points of any  $\phi$-optimal design
 are contained in the set  $\{(1, 0), (0, 1), (1, 1)\}$, and the corresponding weights can now be found by a straightforward calculation.

 In fact
$\phi_{p}$--optimal designs were determined by   \cite{Puke1993}, Section 8.6, who showed
that the $\phi_p$-optimal design for $p\in[-\infty, 1)$ is given by
\[ \xi^*_p=\left\{
\begin{array}{ccc} (1,1) & (1,0) & (0,1) \\
         w(p) & (1-w(p))/2 & (1-w(p))/2
\end{array} \right\}, \]
where  $w(p)=1- {4}/({3+3^{1/(1-p)}})$  if $ p>-\infty$, and   $ w(-\infty)=0$.

For example, the $D$-optimal design $\xi^*_D$, i.e., the $\phi_p$-optimal design with $p=0$,  has masses $1/3, 1/3$ and $1/3$ at the points $(1,1)$, $(1,0)$ and $(0,1)$.
The information matrix of $\xi^*_D$ is given by
\[\M(\xi^*_D)=\left(\begin{matrix}
  2/3&1/3\\
  1/3&2/3
\end{matrix}\right).\]
%From Lemma 6.16 in Pukelsheim (1993), we have the optimal solution of the dual problem for $p=0$ is
%\[\N^*_D=\left(\begin{matrix}
%  1&-1/2\\
%  -1/2&1
%\end{matrix}\right)
%=\left(\begin{matrix}
%  \sqrt{2}/2&\sqrt{2}/2\\
%  -\sqrt{2}/2&\sqrt{2}/2
%\end{matrix}\right)
%\left(\begin{matrix}
%  3/2&0\\
%  0&1/2
%\end{matrix}\right)
%\left(\begin{matrix}
%  \sqrt{2}/2&-\sqrt{2}/2\\
%  \sqrt{2}/2&\sqrt{2}/2
% \end{matrix}\right).\]
The corresponding polytope is obtained for the choice $t=\pi/4$ and  given by
\begin{eqnarray}\label{polytope Ex1}
{\cal P}_{\Z^{*}}=
\Big\{(\lambda_1, \lambda_2): \lambda_1\geq0,\;\lambda_2\geq0, \; \frac{1}{2}\lambda_1+\frac{1}{2}\lambda_2\leq 1,\; 2\lambda_2\leq 1 \Big\}~.
 \end{eqnarray}
% and is displayed in Figure \ref{fig1}.
%The optimal $\bdm{\lambda}^*_D= (3/2,1/2)$ is the point of intersection of the straight line segments $\ell_1: \frac{1}{2}\lambda_1+\frac{1}{2}\lambda_2= 1$  and $\ell_2: 2\lambda_1= 1$.
The two support points $(1,0)$ and $(0, 1)$ correspond to   the same hyperplane defined by the equation   $ \frac{1}{2}\lambda_1+\frac{1}{2}\lambda_2= 1$ since the equalities
 $$
 \big (\cos (\pi/4) x_1 - \sin (\pi/4) x_2  \big)^2=\frac{1}{2} \mbox{~~ and~~}
 \big (\sin (\pi/4) x_1+\cos (\pi/4) x_2  \big )^2=\frac{1}{2}
$$
hold for $(x_1, x_2)=(1,0)$ and $(x_1, x_2)=(0,1)$. The third support point $(1,1)$  corresponds to  the other hyperplane  $2\lambda_1= 1$ because we have
$$
 \big (\cos (\pi/4) x_1 - \sin (\pi/4) x_2  \big)^2 =0 ~~\mbox{ and}  ~~
 \big (\sin (\pi/4) x_1+\cos (\pi/4) x_2  \big )^2=2
 $$
 for $(x_1, x_2)=(1,1).$ Similarly, we consider the $E$-optimal design $\xi^*_E$, i.e., the $\phi_p$-optimal design with $p=-\infty$
 has equal masses at the points $(0,1)$ and $(1,0)$ corresponding to the same  hyperplane
$\frac{1}{2}\lambda_1+\frac{1}{2}\lambda_2= 1$ of the polytope \eqref{polytope Ex1}
and the information matrix is given by
\[\M(\xi^*_E)=\left(\begin{matrix}
  1/2&0\\
  0&1/2
\end{matrix}\right) ~.
\]
%\begin{figure}[t!] % Fig 1
%\centering
%\includegraphics[height=2.5in ,width=3.2in]{DLG_opt.eps}
%\vskip3mm
%\caption{\label{fig1} {  \it  The polytope in (\ref{polytope Ex1}) corresponding to ${N}^*_D$ and ${N}^*_E$. } }
%\end{figure}
}
\end{example}
\vskip3mm

As a direct application of Theorem \ref{thm1}, we obtain sufficient conditions
 for the occurrence of the de la Garza phenomenon in the linear regression model (\ref{LM}), and
 the following corollary gives a bound of the number of support points of an optimal design.

\begin{theorem}\label{Th.4}
Either one of the following conditions is sufficient for the existence of  a $\phi$-optimal design with $k$ support points in the regression model (\ref{LM}).
\begin{description}
\item{(a)} There are no different support points of a design corresponding to the same supporting hyperplane of the polytope as defined in (\ref{hd4}).
\item{(b)} There are no vectors of the same length in the induced design space
$$\mathcal{F}=\{\f(\x) : \x\in\mathcal{X}\}.
$$
\end{description}
\end{theorem}

\begin{corollary}\label{Co1}
If there are at most $N$ vectors of the same length in the induced design space $\mathcal{F}$, then
there exists a $\phi$-optimal design for the regression model (\ref{LM}) with  at most $N k$ support points.
\end{corollary}

\begin{example}\label{wpoly}
{\rm

\noindent
As a first  application  of Theorem \ref{Th.4}, we  consider the de la Garza phenomenon for the weighted   heteroscedastic polynomial regression model
on a compact interval, say $[0,1]$. Optimal design problems in this model have  found considerable attention in the literature \citep[see, for example,][among others]{Fang,Chang2,Dette_Haines2,changlin2006,Chang3,sekido2012}.
To be precise let ${\cal X} = [0,1]$ (or any other compact interval on the non-negative line) and consider the vector of regression functions
\begin{eqnarray}\label{PM}
\f (x)\btheta = \sqrt{\lambda (x)}  (1, x, \ldots ,  x^d)^\top,
\end{eqnarray}
where $\lambda$ is  a positive function on the interval $[0,1]$, which is called  {\it  efficiency function}  in the literature.
It is well known that the information matrix  corresponding to this vector of regression functions is proportional
to the information matrix  in a heteroscedastic polynomial regression model on the interval $[0,1]$, that is
$$
  \mathbb{E} [y(x)] = \theta_0+\theta_1 x+ \cdots+\theta_d x^d ~,~~\mbox{Var}  (y(x)) ={\sigma^{2} \over \lambda (x) }
$$
\citep[see][]{fedorov1972}. If the function $x \to \lambda (x) \|\f(x)\|^2=1+x^2+\cdots+x^{2d}$ is injective    on the   interval $[0,1] $, it follows  from Theorem \ref{Th.4}
that  there exists a $\phi$-optimal design $\xi^*\in\Xi$ supported at  at most $d+1$ points. This
situations occurs in particular, if the function $\lambda$  is increasing,  because the function
$ x \to \|\f(x)\|^2=1+x^2+\cdots+x^{2d}$ is a  strictly increasing function on the interval $[0,1]$.
For the special case $\lambda (x)  \equiv 1$ we obtain  the celebrated de la Garza phenomenon
 \citep[see][]{Garza1954}.
}
\end{example}

\section{Admissibility}
 \label{sec3}
\def\theequation{3.\arabic{equation}}
\setcounter{equation}{0}

In this section, we study the relation between admissibility of a design $\xi$
in the   model (\ref{LM}) and the admissibility of a corresponding  ``conditional  design'' of $\xi$  in a  ``conditional model'' of (\ref{LM}), which will be defined below.
Throughout this paper we call a design $\xi_1$ admissible if there {\bf does not } exist any design $\xi_2$ such that $\M(\xi_1)\neq \M(\xi_2)$ and $\M(\xi_2)\geq \M(\xi_1)$, that is the matrix
 $\M(\xi_2) - \M(\xi_1)$ is nonnegative definite.
%  \HD{A  subclass   of designs  $\Xi_{c} \subset \Xi$ is called  a  {\it complete class},  if for any design  $\xi \in \Xi$    exists a design   $\xi_{c} \in \Xi_{c}$ such that $\M (\xi_{c} ) \geq  \M (\xi )$.}
  For the sake of simplicity, all results  in this section are presented for models with a two-dimensional predictor, but the generalization to the $q$-dimensional
 case with $q \geq 3$ is  straightforward with some additional notation.

% \subsection{Conditional model and marginal model}

To be precise, consider the linear model  (\ref{LM}) with a two-dimensional predictor $\x=(x_1, x_2)$ and define the function
\begin{eqnarray}\label{TF M}
  \mu  (\x ) =   \mu (x_1, x_2)=\sum_{j=1}^pf_j(x_1, x_2)\theta_j=\f^\top(\x)\btheta,
\end{eqnarray}
as the expected response at experimental condition $\x=(x_1, x_2) \in\mathcal{X}$ .
Let $t : \mathcal{X} \to \mathbb{R} $ denote a real-valued function on ${\cal X}$ with range  $\mathcal{T}=\{t(\x) : \x\in \mathcal{X}$\}. The {\it conditional model given $t$}  is
defined on the design space  $\widetilde{\mathcal{X}}(t)=\{ \x\in\mathcal{X}  :  t(\x)=t\}$  (the preimage  of the the set $\{t\}$) and given by
\begin{eqnarray}\label{CM1}
  \widetilde{\mu}_{t}(\x)=\sum_{j=1}^{p_t}\widetilde{f}_{jt}{(\x)}\theta_{jt}
  =\widetilde{\f}_{t}^\top(\x)\widetilde{\btheta}_{t}, \quad \x\in\widetilde{\mathcal{X}}(t),
\end{eqnarray}
where  $\big\{\widetilde{f}_{1t}{(\x)}, \widetilde{f}_{2t}{(\x)},\ldots, \widetilde{f}_{p_t,t}{(\x)}\big\}$ is a set of linearly independent regression functions which spans $\{f_1{(\x)}, f_2{(\x)}, \ldots,f_p({\x})\}$ under the condition that $t(\x)=t$, and $\widetilde{\btheta}_t$ is a $p_t\times 1$  vector of parameters which may depend on $t$.

In the following we are particularly interested in two cases corresponding to the projections on the margins. To be precise
assume that ${\cal X} \subset \mathbb{R}^{2}$ and define  $t_{1}(\x)=x_1$, then for fixed $x_{1}$  the set
$\widetilde{\mathcal{X}}(x_{1}) $  can be identified with  the set ${\mathcal{X}}_2:=\{x_2 : (x_1, x_2)\in\mathcal{X}\}$
and  we obtain  the {\it  conditional model for the second factor $x_2$}  on ${\mathcal{X}}_2$.
Moreover, if  the vector  $\widetilde{\f}_{x_1}(\x)$   in the conditional model \eqref{CM1}
is  independent of $x_1$,
we use the notation $\widetilde{\f}_2(x_2):= \widetilde{\f}_{x_1}(\x)$ and  the resulting model
\begin{eqnarray}\label{Marginal M}
  \widetilde{\mu}_2(x_2)=\widetilde{\f}_2^\top(x_2)\widetilde{\btheta}_2, \quad x_2\in\mathcal{X}_2,
\end{eqnarray}
is called {\it  the marginal model for the second factor}. One can similarly define the conditional model and the marginal model for the first factor.

\begin{example} \label{ex2}  {\rm
To illustrate these ideas we consider the linear model
\begin{eqnarray}\label{IusEx}
  \mu(x_1, x_2)=\theta_0+\theta_1x_1+ \theta_2x_2+\theta_3x_1x_2
\end{eqnarray}
on the design space $\mathcal{X}=[0, 1]^2$.  Consider  the mapping $t(\x)=x_1+x_2$ from the square $[0, 1]^2$ onto the interval $[0,2]$.
For every $t\in[0, 2]$ there are only three independent components among the regression functions $\{1, x_1, x_2, x_1x_2\}$ because of the constraint $x_1+x_2=t$.
Replacing $x_2$ with $t-x_1$ the conditional model can be expressed in the form  (\ref{CM1}) with $\widetilde{\f}_{t}^\top(\x)=(1, x_1, x_1^2)^\top$, where
the conditional design space  $\widetilde{\mathcal{X}}(t)$ can be identified with the interval $[0, t]$ if $t\in(0, 1]$  and
with the interval  $\widetilde{\mathcal{X}}(t)=[t-1, 1]$ if $t\in(1, 2]$. \\
Moreover, the marginal model for the $i$-th factor corresponds to the vector of regression functions ${\widetilde{\f}}_i(x_i)=(1, x_i)^\top$ and the marginal design space
is given by   $\mathcal{X}_i=[0, 1], i=1, 2.$}
\end{example}

\bigskip

For a  design $\xi$ on  the design space  $\mathcal{X}$ we define
$$
\xi_t(t)=\int_{\widetilde{\mathcal{X}}(t)}\xi(d\x)
$$
as the marginal design of $\xi$  on the design space $\mathcal{T}$, then, if  $\xi_t(t) >0 $,  the design $\xi$
 induces  a conditional design $\xi_{\x|t}$  on the design region $\widetilde{\mathcal{X}}(t)$ of the conditional model, which is
 defined by
$$
\xi_{\x|t}(\cdot)=\frac{1}{\xi_t(t)}\xi(\x)~.
$$
  In addition, we define
   $$
   \M_{t}(\xi_{\x|t})=\int_{\widetilde{\mathcal{X}}(t)} \widetilde{\f}_t(\x)\widetilde{\f}_t^\top(\x)\xi_{\x|t}(d\x)
   $$
as  the information matrix of the design $\xi_{\x|t}$  in  the conditional model (\ref{CM1}) and denote  by $\Xi_t$ the set of all approximate designs on the design space
    $\widetilde{\mathcal{X}}(t)$. The following result is proved in the Appendix.

\begin{theorem}\label{Th.6} % Th 6
A necessary condition for the admissibility of a design  $\xi \in \Xi$ in the class   $\Xi$ for the regression model \eqref{TF M} is that   the conditional design $\xi_{\x|t}$ induced by $\xi$
is admissible in the class $\Xi_{t}$    in  the conditional model (\ref{CM1})
 for every $t\in\mathcal{T} $ with $\xi_t(t) >0 $.
\end{theorem}

Furthermore, the following theorem gives a complete subclass and a bound of the number of support points of an optimal design.

\begin{corollary}\label{Th.7}
Assume that the design region is of the form $\mathcal{X}=\mathcal{X}_1\times\mathcal{X}_2$  for some  sets
 $ \mathcal{X}_1,\mathcal{X}_2 \subset \mathbb{R} $ and  suppose  that the marginal models exist for both factors.
Define $ \Xi^A_i$  as  the class of admissible designs for the the $i$-th marginal model ($ i=1, 2$) and denote by
 $\Xi^C$  the   subclass of designs  on ${\cal X} $, where the $i$-th marginal design belong to $\Xi^A_i$ ($i=1, 2$).
 Then the class of all admissible designs    for the model (\ref{TF M}) is a subset of  $\Xi^C$.
 \\
   Moreover, if the admissible designs in $\Xi^A_i$ are based on at most $p_i$ points, $i=1, 2$, then the designs in $\Xi^C$ are based on at most $p_1p_2$ points.

\end{corollary}

%Consequently, we obtain lower bounds for the number of support points of optimal designs for the Kronecker product models and the additive models, respectively.

%\begin{corollary}
%If the admissible designs in $\Xi_l$ for the $l$-th marginal model are based on at most $p_l$ points, $l=1, 2$, then we have the following claims.
%\begin{itemize}
%  \item[(a)]
%  For the Kronecker product type model (\ref{CKP M}),  the designs based on at most $p_1p_2$ points form a complete subclass in $\Xi$;
%  \item[(b)]
%  For the additive model (\ref{add M}), the designs based on at most $\frac{(p_1+p_2+1)^2}{4}$ points form a complete subclass in $\Xi$.
%\end{itemize}
%\end{corollary}

\section{Some applications}
 \label{sec4}
\def\theequation{4.\arabic{equation}}
\setcounter{equation}{0}

In this section we illustrate several  applications of the results in Section \ref{sec2} and \ref{sec3} in the determination of
 locally optimal designs for nonlinear models with a multivariate predictor. To be precise
we consider the common nonlinear regression models with $q$ factors
\begin{eqnarray}\label{NLModel}
  \mathbb{E} [ y(\x) ] =\eta(\x,\btheta), \quad\x\in\mathcal{X}\subset\mathbb{R}^q,
\end{eqnarray}
where $y(\x)$ is a normal distributed random variable   with constant variance, say $\sigma^{2} >0$
and observations at different experimental conditions are assumed to be independent.
We further assume that the (non-linear) regression function $\eta(\x,\btheta)$ is continuously differentiable with respect to the parameter $ \btheta $ and
define
\begin{eqnarray} \label{hd11}
  \f({\x},\btheta)=\nabla\bdm{\eta}(\x,\btheta)=\left(\frac{\partial\bdm{\eta}(\x,\btheta)}{\partial\theta_1}, \ldots, \frac{\partial\bdm{\eta}(\x,\btheta)}{\partial\theta_k}\right)^\top,
\end{eqnarray}
as the gradient of $\eta$ with respect to the parameter $ \btheta $. The information matrix of a design $\xi$ for model (\ref{NLModel}) is given by
\begin{equation}\label{Informatix Matrix for NLM}
{\M}(\xi, \btheta)=\dint_{\mathcal{X}}\f({\x},\btheta){\f}^\top({\x},\btheta) \xi({d\x}).
\end{equation}
If $n$ observations are taken according to an approximate design (applying an appropriate rounding procedure)
 it is well know, that under standard assumptions,  the covariance matrix  of the  maximum likelihood estimate of the parameter $ \btheta $ is
approximately given by   the matrix  $\sigma^{2}/n {\M}^{-1}(\xi, \btheta)$ and a locally  optimal design maximizes an information  function
of the matrix ${\M}(\xi, \btheta)$.  Consequently, the results of the previous sections can be used to characterize properties of admissible designs
for locally optimal design problems, where the vector of regression function is given by the gradient
$\f (\x, \btheta )$ defined in  \eqref{hd11}. We illustrate this in a few examples.

\subsection{Exponential regression}
\label{ex5}
 \rm
\cite{detmelpep2006,DMW2006}  studied optimal designs for the exponential regression model
\begin{eqnarray}\label{ERM}
  \bdm{\eta}(x,\btheta)=\sum_{l=1}^La_l\exp(-\lambda_lx),\quad x\in \mathcal{X}=[0, \infty),
\end{eqnarray}
where the vector of parameters is given by $\btheta=(a_1, \ldots, a_L, \lambda_1, \ldots, \lambda_L)^\top$ with $a_l\neq0$, $l=1,\ldots,L,$ and $0<\lambda_1<\cdots<\lambda_L.$
For $L = 1$ or 2,  \cite{DMW2006}  showed that there exists  a locally $D$-optimal design based on $2L$ points.
Moreover,  for $L \geq 3$ they defined   $\widetilde{\bdm{\lambda}}=(\widetilde{\lambda}_1, \ldots,\widetilde{\lambda}_L)^\top$  as the  vector
with components satisfying $0<\widetilde{\lambda}_1<\cdots<\widetilde{\lambda}_L$ and $\widetilde{\lambda}_{i+1}=(\widetilde{\lambda}_{i}+\widetilde{\lambda}_{i+2})/2$, $i=1,\ldots, L-2$, and  showed  that for  any vector $\widetilde{\bdm{\lambda}}$ of this type  the  existence of  a neighbourhood $\mathcal{U}$ of $\widetilde{\bdm{\lambda}}$, such that for all vectors $\bdm{\lambda}
=(\lambda_{1}, \ldots , \lambda_{L})^{\top}\in\mathcal{U}$, there
exists a locally $D$-optimal design for the parameter $\btheta$  which is supported on $2L$  points. Moreover, they pointed out that numerical results indicate that the set of parameter vectors $\bdm{\lambda}$ for which the locally D-optimal design is minimally supported is usually very large.
\cite{YS2012}  established similar conclusions for optimal designs with respect to other   criteria in  the   cases  $L=2$  (here the condition $\lambda_1/\lambda_2<61.98$
is sufficient for the existence of locally optimal design supported at $4$ points) and $L=3$ (here the conditions   $2\lambda_2=\lambda_1+\lambda_3$ and $\lambda_2/\lambda_1<23.72$
imply the existence of an optimal design supported at $6$ points).

We now  extend these results in a non-trivial manner using the methodology developed in Section \ref{sec2}.

\begin{theorem}  \label{corerm}
Consider the exponential regression model \eqref{ERM} on the interval $[0,\infty)$, where
$0 < \lambda_1 < \lambda_2 < \ldots < \lambda_L$  and $a_l\neq0$ $l=1,\ldots,L$.
If  the parameters satisfy $\lambda_i \geq \frac {|a_{i}|}{2}$  for all $i=1, \ldots , L$, any $\phi$-optimal design is supported at at most $2L$ points.
%  In particular, all $\phi_p$-optimal designs for $p \in [-\infty,0]$ for the model \eqref{ERM} are supported at $2L$ points. \\
Moreover, for  all  $0 < \lambda_1 < \lambda_2 < \ldots < \lambda_L$  and $a_1, \ldots , a_{L}\neq0$
there  exists a locally $D$-optimal design supported at $2L$ points.
\end{theorem}

\subsection{Exponential regression models with two factors}

\cite{ROM2015}  considered the maximin  optimal design problem for  the  two-factor exponential growth model
\begin{eqnarray}\label{EGM}
  \bdm{\eta}(\x,\btheta)=\theta_0+\exp(-\theta_1x_1)+\exp(-\theta_2x_2),
\end{eqnarray}
($\theta_j\geq1$, $j=1,2$) on the square $\mathcal{X}=[0, 1]^2$, which have numerous applications in
biological and agricultural sciences. In  this model the  gradient  of  the function $\bdm{\eta}(\x,\btheta)$ in (\ref{EGM})  is given by
\begin{eqnarray}\label{RF_NE2}
  \f(\x,\btheta)=\left(1, -x_1\exp(-\theta_1x_1),-x_2\exp(-\theta_2x_2)\right)^\top
\end{eqnarray}
and the two vectors of regression functions corresponding to the marginal models of (\ref{RF_NE2}) are obtained as
\begin{eqnarray}\label{EGM_M1}
  \widetilde{\f}_i(x_i, \btheta)=\left(1, x_i\exp(-\theta_ix_i)\right)^\top, \quad  i=1,2.
\end{eqnarray}
The admissible designs for the marginal models  are supported at  the points  $\{0, 1/\theta_i\}$  ($i=1,2$), and
it now follows from  Theorem \ref{Th.6} and  Corollary  \ref{Th.7} that  the admissible  designs   for  model (\ref{EGM}) are contained in the class of all  designs supported
at most $4$  points from  the set  $\{(0, 0), (0, 1/\theta_2),$ $(1/\theta_1, 0), (1/\theta_1, 1/\theta_2)\}.$
For example, a straightforward optimization shows that
the locally  $D$-optimal design puts masses $1/4$ at all four points.

\vskip3mm

Similarly, admissible designs can be determined for   the  two-factor exponential model
\begin{eqnarray}\label{TDEM}
  \bdm{\eta}(\x,\btheta)=\theta_0\exp{(\theta_1x_1+\theta_2x_2)},
\end{eqnarray}
where $\theta_j>0$, $j=1,2,3$ and the design space is given by  $\mathcal{X}=[0, b_1]\times[0, b_2].$
\cite{GMS2018}  investigated the locally $D$-optimal designs for this model by means of a general equivalence theorem and showed that  locally $D$-optimal designs
are   supported at at most $4$  points.

The gradient of the function $\bdm{\eta}(\x,\btheta)$ in  model (\ref{TDEM})  is given by
\begin{eqnarray}
\label{RF_NE3}
  \f (\x,\btheta)= \exp{(\theta_1x_1+\theta_2x_2)} \big ( 1  , \theta_0 x_1,\theta_0x_2 \big )^\top.
\end{eqnarray}
For given $\btheta$ let $t=\theta_1x_1+\theta_2x_2$  and define the matrix
\[ %\begin{eqnarray}
\C(t, \btheta)=\begin{pmatrix}
  \exp(t)&0&0\\
  0&\theta_0\exp{(t)}&0\\
  0&0&\theta_0\exp{(t)}
\end{pmatrix},
%\end{eqnarray}
\]
then the vector of  regression functions corresponding to the conditional model of (\ref{RF_NE3}) is given by
\begin{eqnarray}\label{TDEM_CM}
  \widetilde{\f}_t(\x)=\left(1, x_1, x_2\right)^\top
\end{eqnarray}
and the design space for   the conditional model is given by $\widetilde{\mathcal{X}}(t)=\{\x\in\mathcal{X}  : \theta_1x_1+\theta_2x_2=t \}$.
For every $t\in\mathcal{T}=\{t=\theta_1x_1+\theta_2x_2 : \x\in\mathcal{X}\}$, it is easy to see that the admissible design for the conditional model (\ref{TDEM_CM}) on the design region $\widetilde{\mathcal{X}}(t)$ is supported at the two end points of the line segment $\ell: \theta_1x_1+\theta_2x_2=t$, $\x\in\mathcal{X}$. Therefore, by Theorem
\ref{Th.6}, the admissible designs for the model (\ref{TDEM}) are supported on the boundary of the design region $\mathcal{X}$.

Moreover, the two  marginal models  of (\ref{RF_NE3}) exist with corresponding  vectors of regression functions
 given by
\begin{eqnarray}\label{TDEM_M1}
  \widetilde{\f}_i(x_i, \btheta)=\left(\exp(\theta_ix_i), x_i\exp(\theta_ix_i)\right)^\top, \quad  i=1,2.
\end{eqnarray}
The admissible designs for the $i$-th marginal model are supported at two points, one of which is $b_i$  ($i=1,2$)  \citep[see][Theorem 5]{YS2009}. It
 now follows from   Theorem \ref{Th.6} and
 Corollary  \ref{Th.7} that  the admissible  designs   for  model (\ref{TDEM}) are contained in the class of all  designs supported at the
$ 4$  points, $(b_1, b_2)$, $(b_1, x_2)$, $(x_1, b_2)$, $(x_1, x_2)$, where $x_i\in [0, b_i)$ and  the point $(x_1, x_2)$ is a  boundary point of
$\mathcal{X}=[0, b_1]\times[0, b_2].$

\subsection{Mixture of exponentials and polynomials}

\cite{ROM2015}  considered the maximin  optimal design problem for  the  model
\begin{eqnarray}\label{MEPM}
  \bdm{\eta}(\x,\btheta)=\theta_0+\theta_1 x_1+\theta_2 x_1^2+\exp(-\theta_3x_2),
  \end{eqnarray}
  which was used in \cite{Langseth2012} for approximating the potentials associated with general hybrid Bayesian networks.
The parameter $\theta_3$ is assumed to be positive,  design space is given by $\mathcal{X}=[-1, 1]\times[0, 2]$ and the
 gradient  of  the function $\bdm{\eta}(\x,\btheta)$ in (\ref{MEPM})  is obtained as
\begin{eqnarray}\label{RF_NE4}
  \f(\x,\btheta)=\left(1, x_1, x_1^3,-x_2\exp(-\theta_3x_2)\right)^\top~.
\end{eqnarray}
The marginal  models of (\ref{RF_NE4}) exist with corresponding  vectors of regression functions
given by
\begin{eqnarray}\label{MEPM_M1}
  \widetilde{\f}_1(x_1, \btheta)&=&(1, x_1, x_1^3)^\top, \quad\mathcal{X}_1=[-1, 1]
\\
\label{MEPM_M2}
  \widetilde{\f}_2(x_2, \btheta)&=&(1, x_2\exp(-\theta_3x_2))^\top, \quad \mathcal{X}_2=[0, 2].
\end{eqnarray}
The admissible designs for the marginal model (\ref{MEPM_M1}) are supported at at most $4$ points including end points $-1$ and $1$
\citep[see][Theorem 8]{Yang2010}. In addition, it follows from Corollary \ref{Co1} that the other two support points, say $u^*$ and $v^*$, satisfy the condition $ \|\widetilde{\f}_1(u, \btheta)\|=\|\widetilde{\f}_1(v, \btheta)\|$, which implies $u^*=-v^*$.
For the marginal model (\ref{MEPM_M2}) the admissible designs are supported at the points  $\{0,x_2^*\}$,  where $x_2^*=\min\{1/\theta_3, 2\}$. %and $x_2^*=2$ if $\theta_3\leq0$.
It now follows from  Corollary  \ref{Th.7} that  the admissible  designs   for  model (\ref{MEPM}) are contained in the class of all  designs supported at the
$ 8$  points, $(\pm 1, 0)$, $(\pm 1, x_2^*)$, $(\pm u^*, 0)$, $(\pm u^*, x_2^*)$, where $u^*\in [0, 1)$. For example, the locally $D$-optimal design for model (\ref{MEPM}) is equally supported at $(\pm 1, 0), (\pm \sqrt{3}, 0), (\pm 1, x_2^*), (\pm \sqrt{3}, x_2^*)$.

\section{Appendix: proofs}
 \label{sec5}
\def\theequation{5.\arabic{equation}}
\setcounter{equation}{0}

\subsection{Preliminaries}
 \label{sec51}
In this section we recall some general results from optimal design theory which will be used in subsequent proofs. For more details the reader is referred to the monograph of
\cite{Puke1993}.

The  polar function $\phi^{\infty}:$ NND($s$)$\rightarrow[0;\infty)$
 of an information function  $\phi$: PD($s$)$\rightarrow (0,\infty)$  is defined by
  \begin{equation}
    \phi^{\infty}(\D)=\inf_{\C \in {\rm PD}(s) }\frac{\trace(\C\D)}{\phi(\C)} ~.
  \end{equation}
{For every information function function $\phi$   the corresponding polar function $\phi^{\infty}$ is isotonic relative to the Loewner ordering.} Define
\begin{eqnarray}\label{Set N}
  \mathcal{N}=\{\N\in \textmd{NND}(k): \f^\top(\x)\N\f(\x)\leq 1 \;\; \forall \x\in\mathcal{X}\},
\end{eqnarray}
then a duality relation    of the optimal design problem can be established [see \cite{Puke1993}, Theorem 7.12], that is
  \begin{eqnarray}
    \max_{\M\in\mathcal{M}(\Xi)}\phi(\M)=\min_{N\in\mathcal{N}}\frac{1}{\phi^{\infty}(\N)}.
  \end{eqnarray}
In particular, an information matrix $\M\in\mathcal{M}(\Xi)$ is optimal for $\btheta$ in $\mathcal{M}(\Xi)$ if and only if there exists a matrix $\N\in\mathcal{N}$ such that
\begin{eqnarray} \label{hd1}
\phi(\M)=\frac{1}{\phi^{\infty}(\N)},
\end{eqnarray}
and  two matrices $\M\in\mathcal{M}(\Xi)$ and $\N\in\mathcal{N}$ satisfy \eqref{hd1}  if and only if  the conditions
\begin{eqnarray}
  \trace(\M\N)&=&1,\\
  \phi(\M){\phi^{\infty}(\N)}&=&\trace(\M\N)
\end{eqnarray}
hold. An application of this result  yields the famous general  equivalence theorem in optimal design theory.
\begin{theorem}[\cite{Puke1993}, Theorem 7.17]\label{Th.2} % Th.2
A positive definite   information matrix  $\M^{*}\in\mathcal{M}(\Xi)$ is  $\phi$-optimal for $\btheta$ in $\mathcal{M}(\Xi)$ if and only if there exists a nonnegative definite
$k\times k$ matrix $\N \in {\cal N} $ that solves the polarity equation
\begin{eqnarray*}
  \phi(\M^{*}){\phi^{\infty}(\N)}=\trace(\M^{*}\N)=1
\end{eqnarray*}
and that satisfies the normality inequality
\begin{eqnarray}  \label{HD2}
  \f^\top(\x)\N\f(\x)\leq 1 \;\; \forall \x\in\mathcal{X}.
\end{eqnarray}
Moreover, if $\M^{*} $ is optimal for $\btheta$ in $\Xi$, there is equality for  any support point $\x_i$ of any $\phi$-optimal design $\xi\in\Xi$, that is any design with  $\M^{*}  = \M (\xi)$.
\end{theorem}

\subsection{Proof of Theorem \ref{thm1}}
 \label{sec52}

The matrix   $\N$ in NND($k$) has an eigenvalue decomposition
  \begin{eqnarray*}
    \N=\Z_\N\bLam_\N \Z_\N^\top,
  \end{eqnarray*}
where $\bLam_\N=$ diag $(\lambda_{\N 1}, \ldots, \lambda_{\N k})$ is a diagonal matrix,   the eigenvalues $\lambda_{N1}, \ldots, \lambda_{Nk}$ of $\N$ are counted with their respective multiplicities, and $\Z_\N=(\z_{\N 1}, \dots, \z_{\N k})$ is an orthogonal matrix with eigenvectors corresponding to the eigenvalues.
Denote by $\mathcal{S}_\Z$ the subset of NND($k$) consisting of matrices which permit eigenvalue decomposition with the same orthogonal matrix $\Z$, i.e.,
$$
\mathcal{S}_\Z=\{\N : \N=\Z\bLam \Z^\top, \bLam=\textmd{diag}(\lambda_1, \ldots, \lambda_k), \lambda_i\geq0, i=1, \ldots, k\}.
$$
 Then we can express NND($k$) as
 $$
 \text{NND}(k)=
 \bigcup_{\Z\in O(k)}\mathcal{S}_\Z,
 $$
  where $O(k)$ is the set of all $k \times k$ orthogonal matrices.

Furthermore, let
\[\h_\Z(\x)=\Z^\top\f(\x)=(\z_1^\top\f(\x), \ldots, \z_k^\top\f(\x))^\top=(h_{\Z 1}(\x), \ldots, h_{\Z k}(\x))^\top, \]
 then the set $\mathcal{N}$ in (\ref{Set N}) can be represented as
\[\begin{array}{lll}
\mathcal{N}
  &=&\bigcup_{\Z\in O(k)}\{\N\in \mathcal{S}_\Z: \f^\top(\x)\N\f(\x)\leq 1 \;\; \forall \x\in\mathcal{X}\} \crno
   &=&\bigcup_{\Z\in O(k)}\{\N\in \mathcal{S}_\Z: \h_\Z^\top(\x)\bLam\h_\Z(\x)\leq 1 \;\; \forall \x\in\mathcal{X}\} \crno
   &=&\bigcup_{\Z\in O(k)}\{\N\in \mathcal{S}_\Z: \h_{\Z 1}^2(\x)\lambda_1+\cdots+\h_{\Z k}^2(\x)\lambda_k\leq 1 \;\; \forall \x\in\mathcal{X}\} \crno
   &\widehat{=}&\bigcup_{\Z\in O(k)}\mathcal{N}_\Z ,
 \end{array}
\]
where the last line defines the set  $\mathcal{N}_\Z $ in an obvious manner.
The optimal solution of the dual problem must occur on some subset, say $\mathcal{N}_{\Z^*}$. Moreover, the dual problem on any subset $\mathcal{N}_\Z$ can be viewed as an extremum problem of a multivariate function defined on the convex polytope   \eqref{hd4}, that is
$$
\{(\lambda_1, \ldots, \lambda_k): \lambda_i\geq0, \; i=1, \ldots, k, \;\h_{\Z 1}^2(\x)\lambda_1+\cdots+\h_{\Z k}^2(\x)\lambda_k\leq 1 \;\; \forall \x\in\mathcal{X}\}.
$$
Note that the polar function $\phi^\infty$ is isotonic, hence the optimal solution must be attained at a boundary point of the polytope. Let $\bdm{\lambda}^*=(\lambda_1^*, \ldots, \lambda_k^*)$ be the vector corresponding to the extremum,   then there are at most $k$ effective constraints, say
 $$
 \ell_i: c_{i1}\lambda_1+\cdots+c_{ik}\lambda_k = 1, \quad i=1, \ldots,  k,
  $$
  which must be satisfied by  $\bdm{\lambda}^*$. Corresponding to $\bdm{\lambda}^*$ we define $\N^*=\Z^*\bLam^*\Z^{*\top}$ with $\bLam^*=$ diag $(\lambda_1^*, \ldots, \lambda_k^*)$,
 then $\N^*$ is the optimal matrix of the dual problem.

By Theorem \ref{Th.2}, we have
 \[1=\f^\top(\x^*)\N^*\f(\x^*)=h_{\Z^*1}^2(\x^*)\lambda_1^*+\cdots+h_{\Z^*k}^2(\x^*)\lambda_k^*\]
for any support point $\x^*$ of a $\phi$-optimal design $\xi^*$, which implies that \[(h_{\Z^*1}^2(\x^*), \ldots, h_{\Z^*k}^2(\x^*))\in\{\c_i=(c_{i1}, \ldots, c_{ik})^\top, i=1, \dots, k\},\]
and $h_{\Z^*1}^2(\x^*)\lambda_1+\cdots+h_{\Z^*k}^2(\x^*)\lambda_k= 1$ is a supporting hyperplane of the polytope (\ref{hd4}) with $\Z=\Z^*$.
Consequently, the support points of $\xi^*$ can
be divided into $k$ sets, say $\{\x_{ij},   j=1, \ldots, m_i\}$, according to the vectors $\c_1, \ldots, \c_k$. Moreover, the support points in the $i$-th set satisfy
\[  %\begin{eqnarray}\label{Key Eq}
 (h_{\Z^*1}^2(\x_{ij}), \ldots, h_{\Z^*k}^2(\x_{ij}))=\c_i, \quad j=1, \ldots, m_i,
 %\end{eqnarray}
\]
which yields
 \[(\z_l^{*\top}\f(\x_{ij}))^2=c_{il}, \quad l=1, \ldots, k; \; j=1, \ldots, m_i.\]
Therefore, the vectors $\f(\x_{ij}), j=1, \ldots, m_i$ share the same length since
\[\|\f(\x_{ij})\|^2=\|\Z^{*\top}\f(\x_{ij})\|^2
=\sum_{l=1}^k(\z_l^{*\top}\f(\x_{ij}))^2=\|\c_i\|^2,\]
which completes the proof of Theorem \ref{thm1}.

\subsection{Proof of Theorem \ref{Th.6}}
For fixed $t$,   there exists a full column-rank matrix, say $\C(t)$, such that
\[ %\begin{eqnarray}
 \f(\x)=\C(t) \widetilde{\f}_t(\x)
%\end{eqnarray}
\]
on the design space $\widetilde{\mathcal{X}}(t)$,
since the elements of the vector $\widetilde{\f}_t(\x)$ are   linearly independent.
Suppose there exists some $t_*\in\mathcal{T}$  with $\xi_{{t}}(t_*) >0 $
such that the conditional design $\xi_{\x|t_*}$ is inadmissible for the conditional model (\ref{CM1}). Then there exists a design
$\bar{\xi}_{\x|t_*}$ in a set of all conditional designs $\Xi_{t_*}$ satisfying
\[ %\begin{eqnarray}
  \M_{t_*}(\bar{\xi}_{\x|t_*}){\geq}\M_{t_*}(\xi_{\x|t_*}) ~ \mbox{ and }~
  \M_{t_*}(\bar{\xi}_{\x|t_*}){\neq}\M_{t_*}(\xi_{\x|t_*}).
%\end{eqnarray}
\]
Let $\bar{\xi}$ be the design obtained by replacing  the conditional design $\xi_{\x|t_*}$ of $\xi$ with $\bar{\xi}_{\x|t_*}$, then we have
\begin{eqnarray*}
   \M(\xi)&=&\int_{\mathcal{X}}\f(\x)\f^\top(\x) \xi(d\x)\\
  &=&\int_{\mathcal{T}}\int_{\widetilde{\mathcal{X}}(t)}\f(\x)\f^\top(\x) \xi_{\x|t}(d\x)\xi_t(dt)\\
  &=&\int_{\mathcal{T}}\int_{\widetilde{\mathcal{X}}(t)}\C(t) \widetilde{\f}_t(\x)\widetilde{\f}_t^\top(\x)\C^\top(t) \xi_{\x|t}(d\x)\xi_t(dt)\\
  &=&\int_{\mathcal{T}}\C(t)\left[\int_{\widetilde{\mathcal{X}}(t)}\widetilde{\f}_t(\x)\widetilde{\f}_t^\top(\x) \xi_{\x|t}(d\x)\right]\C^\top(t)\xi_t(dt)\\
   &=&\int_{\mathcal{T}}\C(t)\M_t(\xi_{\x|t})\C^\top(t)\xi_t(dt)\\
   &\leq\atop\neq&\int_{\mathcal{T}}\C(t)\M_t(\bar{\xi}_{\x|t})\C^\top(t)\xi_t(dt)\\
   &=&\int_{\mathcal{T}}\C(t)\M_t(\bar{\xi}_{\x|t})\C^\top(t)\bar{\xi}_t(dt)\\
     &=&\M(\bar{\xi}).
\end{eqnarray*}
Therefore, the design $\xi$ would be  inadmissible in the class $\Xi$ for the model (\ref{TF M}), and this  contradiction completes the proof of Theorem \ref{Th.6}.

 \subsection{Proof of Theorem \ref{corerm}}

The fist part follows directly from  an application of Theorem  \ref{Th.4}(b).  To be precise, note that the gradient  of   the function $\bdm{\eta}(\x,\btheta)$  with respect to the parameter $ \btheta $ in model  (\ref{ERM}) is given by
\begin{eqnarray}\label{RF_NE1}
 \f(x,\btheta)=\left(\exp(-\lambda_1x), -a_1x\exp(-\lambda_1x),\ldots,\exp(-\lambda_Lx), -a_Lx\exp(-\lambda_Lx)\right)^\top
% &=& \diag(1, -a_1, \ldots, 1, -a_L )  \g(x,\btheta) ,
\end{eqnarray}
%where
%$$
% \g(x,\btheta) = (\exp(-\lambda_1x), x\exp(-\lambda_1x),\ldots,  \exp(-\lambda_Lx),x\exp(-\lambda_Lx))^\top.
%$$
It is easy to see that the function $x \to \|\f(x,\btheta)\|^2=\sum_{i=1}^L\exp(-2\lambda_ix)(1+a_{i}^{2}x^2)$ is a strictly decreasing  function on  the interval $[0, \infty)$ if $\lambda_i\geq\frac{|a_{i}|}{2}$  for all $i=1, \ldots , L$, since in this case the derivative of
function $\|\f(x,\btheta)\|^2$ is non-positive.
\\
For the statement regarding the $D$-optimality criterion
recall  that the parameter vector in \eqref{RF_NE1} is given by
 $\btheta=(a_1, \ldots, a_L, \lambda_1, \ldots, \lambda_L)^\top$. Let $c>0$ be any constant  and note that  the vector of regression functions satisfies
\[\begin{array}{lll}  \f(x,\btheta)
 &=& \Q  \g(x,\btheta) ,
\end{array}
\]
where
$$
 \g(x,\btheta) = (\exp(-\lambda_1x), c x\exp(-\lambda_1x),\ldots,  \exp(-\lambda_Lx), c x\exp(-\lambda_Lx))^\top.
$$
and the matrix $\Q$ is given by $\Q=\diag(1, -a_1/c, \ldots, 1, -a_L/c ) $.
Observing the relation
$$
\int_{\cal X}  \f(x,\btheta)  \f^{\top }(x,\btheta)  \xi (d\x)  = \Q \int_{\cal X}  \g(x,\btheta)  \g^{\top }(x,\btheta)  \xi (d\x) \Q^{\top}
$$
it is easy to see that a design is $D$-optimal for the regression model  \eqref{LM} with vector $\f$ if and only if it is $D$-optimal
for the regression model  \eqref{LM} with vector $\g$. However, from the first part of the proof this design is the locally $D$-optimal design
if $\lambda_{i} \geq  c/2$ for all $i=1, \ldots ,  L$. As the constant $c>0$ can be chosen arbitrarily, the assertion of Theorem \ref{corerm}
follows.

\bigskip\bigskip

{\bf Acknowledgements}
Dr. Liu and Dr. Yue were partially supported  by the National Natural Science Foundation of China under Grants 11871143, 11971318.
The work of Dr. Dette has been supported in part by the
Collaborative Research Center ``Statistical modeling of nonlinear
dynamic processes'' (SFB 823, Teilprojekt C2) of the German Research Foundation
(DFG).

\end{document}